\documentclass[12pt,twoside]{article}

\usepackage{amsmath,amsfonts,amsthm,amssymb,amsopn,amstext,amscd}
\usepackage{enumerate}
\usepackage{graphicx}
\usepackage{float}
\usepackage{xcolor}
\usepackage{longtable}
\usepackage[utf8]{inputenc}

\usepackage{setspace} 
\usepackage{enumitem}
\usepackage{hyperref}
\usepackage{dsfont}
\usepackage{url}
\allowdisplaybreaks

\newtheorem{teorema}{Theorem}[section]
\newtheorem{teoremaA}{Theorem A\hspace*{-0.2cm}}

\newtheorem{proposicion}{Proposition}[section]

\newtheorem{lema}{Lemma}[section]
\newtheorem{lemaA}{Lemma A\hspace*{-0.2cm}}
\newtheorem{nota}{Remark}[section]

\newcommand*\samethanks[1][\value{footnote}]{\footnotemark[#1]}

\begin{document}
\pagenumbering{arabic}
\singlespace

\title{Diffusion approximation of controlled branching processes using limit theorems for {random step processes}}
\author{Miguel Gonz\'alez\ \footnote{Department of Mathematics, Faculty of Sciences and Instituto de Computaci\'on Cient\'ifica Avanzada, University of Extremadura, Badajoz, Spain. E-mail address: \url{mvelasco@unex.es}. ORCID: 0000-0001-7481-6561.}
\and Pedro Mart\'in-Ch\'avez\ \samethanks[1]\footnote{ Department of Mathematics, Faculty of Sciences, University of Extremadura, Badajoz, Spain. E-mail address: \url{pedromc@unex.es}. ORCID: 0000-0001-5530-3138.}
\and In\'es del Puerto\samethanks[1]\  \footnote{Department of Mathematics, Faculty of Sciences and Instituto de Computaci\'on Cient\'ifica Avanzada, University of Extremadura, Badajoz, Spain. E-mail address: \url{idelpuerto@unex.es}. ORCID: 0000-0002-1034-2480.} }

\maketitle
\begin{abstract}
A controlled branching process (CBP) is a modification of the standard Bienaymé-Galton-Watson process in which the number of progenitors in each generation is determined by a random mechanism. We consider a CBP starting from a random number of initial individuals. The main aim of this paper is to provide a Feller diffusion approximation for critical CBPs. A similar result by considering a fixed number of initial individuals by using operator semigroup convergence theorems has been previously proved in \cite{sriram}. An alternative proof is now provided making use of limit theorems for {random step processes}.
\end{abstract}

\noindent\textbf{Keywords}: Controlled branching processes; Weak convergence theorem; Martingale Differences; Diffusion processes; Stochastic differential equation; {random step processes}.

\section{Introduction}\label{sec:intro}
Let $\left\{X_{n,j}: n=0,1,\ldots; j=1,2,\ldots\right\}$ be a sequence of
independent and identically distributed (i.i.d.), non-negative and integer-valued random variables defined on a probability space $(\Omega,\cal{F},P)$. Let also $\{\phi_{n}(k): k=0,1,\ldots \}$, for $n=0,1,\ldots$, be a sequence of stochastic processes which consist of independent non-negative
integer-valued random variables {on $(\Omega, \mathcal{F},P)$} with the same one-dimensional distributions. Furthermore, let us assume that $\{X_{n,j}: n=0,1,\ldots; j=1,2,\ldots\}$ and
$\{\phi_{n}(k): n=0,1,\ldots; k=0,1,\ldots\}$ are independent.

A controlled branching process (CBP) is defined recursively as
\begin{equation}{\label{a002}}
 Z_{n} = \sum^{\phi_{n-1}(Z_{n-1})}_{j=1} X_{n-1,j}, \qquad n = 1,2,\ldots,
\end{equation}
where $\sum_{j=1}^0$ is defined as 0 and $Z_{0}$ is a non-negative, integer-valued,
square-integrable random variable which is independent of
$\left\{X_{n,j}: n=0,1,\ldots; j=1,2,\ldots\right\}$ and $\{\phi_{n}(k): n=0,1,\ldots; k=0,1,\ldots\}$.

Here, $Z_{n}$ denotes the size of the $n$-th generation of a
population and $X_{n-1,j}$ is the offspring size of the $j$-th
individual in the $(n-1)$-th generation. We will assume that
the mean $m=E[X_{0,1}]$ and variance $\sigma^{2}=Var[X_{0,1}]$ are both finite.

The class of CBPs is a very general family of stochastic processes that collect as  particular cases  the simplest  branching model, the standard Bienaym\'e–Galton–Watson (BGW) process, by considering $\phi_n(k) = k$ a.s. for each $k$, or {a} branching processes
with immigration, by setting $\phi_n(k)  = k+I_n$, where $\{I_n\}_{n\geq 0}$
are i.i.d. random variables ({writing in this way the immigrants give rise to offspring  at the same generation as their arrival and with the same probability law as $X_{0,1}$)},
among others.  The  monograph \cite{gpy} provides an extensive description of its probabilistic theory.

The research of functional weak limit theorems for branching processes arises
a lot of interest since many years ago. It was firstly formulated for a BGW process by \cite{feller} and proved by \cite{jirina} and \cite{lindvall}. These results have been extended to another classes of branching processes. For instance, a wide literature exists around weak convergence results for  branching processes with immigration (BPI) since the pioneer work by \cite{ww}, see also \cite{bbp21} and references therein. In this paper we focus our attention on a weak convergence theorem for a {critical} CBP with a random initial number of individuals { and assuming finite second order moment on the this initial value}. A similar result was already established for a single CBP in \cite{sriram}, and for an array of CBPs in \cite{b3}, by assuming fixed initial numbers of progenitors  using infinitesimal generators results for their proofs. Inspired in the paper \cite{bbp21} on BPI we will use limit theorems for random step processes
towards a diffusion process provided in \cite{ip10} to obtain an alternative proof. {The scheme of it follows similar steps to the ones in \cite{bbp21}.   An important feature of a CBP is that   the value of $Z_n$ conditioned on the knowledge of the previous generation, $Z_{n-1}=k$,   is a random sum of random variables, namely $\sum_{j=1}^{\phi_n(k)}X_{n-1,j}$, instead of a  non-random sum as  in the case  of a BPI. This leads to handle the proofs of each steps using  conditioning arguments different from those  used in  \cite{bbp21}.  }

Apart from this introduction, the paper is organized as follows. In Section \ref{section2} we provide
the notation and some auxiliary results about the behaviour of the first and second moments of the process.
Section \ref{sec:main} gathers the main theorem. For the ease of reading the paper, additional results  are presented in the Appendix.

\section{Notation and auxiliary results}\label{section2}
We denote, for $ k = 0, 1, \ldots,$
\begin{eqnarray*}
    \varepsilon(k) & = & E[\phi_{n}(k)], \\
    \nu^{2}(k) & = & Var[\phi_{n}(k)],
\end{eqnarray*}
and assume all finite.
It is easy to  obtain that for $n =  1,2, \ldots$,
    \begin{eqnarray}
    E[Z_{n}|{\cal F}_{n-1}] & = & m \varepsilon(Z_{n-1}),\label{esperanza}\\
    Var[Z_{n}|{\cal F}_{n-1}] & = & \sigma^{2} \varepsilon(Z_{n-1}) + m^{2} \nu^{2}(Z_{n-1})\label{varianza},
    \end{eqnarray}
where ${\cal F}_{n}$ is the $\sigma-$algebra generated by the random variables $Z_{0}, Z_{1},\ldots, Z_{n}$, $n\geq 1 $ {(see Proposition 3.5 in \cite{gpy})}.

We introduce the quantities
\begin{equation}{\label{taumk}}
\tau_{m}(k)= E[Z_{n+1}Z^{-1}_{n}|Z_{n} = k] = m
\varepsilon(k)k^{-1}, \ k\geq 1.
\end{equation}

The quantity $\tau_{m}(k)$ represents a mean growth
rate. Intuitively, it can be interpreted as an average offspring
per individual for a generation of size $k$.

Assuming that {$\lim_{k\to\infty}\tau_m(k)=\tau_m$ exists, the process can be classified as:
$$\tau_m <1 \quad\mbox{{subcritical}}; \ \tau_m =1 \quad\mbox{{critical}};\  \tau_m>1 \quad\mbox{{supercritical}}.$$}

We are interested in critical CBPs that satisfy the following hypotheses:
\begin{itemize}
    \item[A1)] $\tau_m(k) = 1 + k^{-1}\alpha$, $\ k >0$, $\alpha>0$,
    \item[A2)] $\nu^2(k) = O(k^{\beta})$, $\beta<1$, as $k\to\infty$.
\end{itemize}

{ The behavior of  critical CBPs was studied in
\cite{b2}.   Assuming that $P(\phi_0(0)=0)=1$, i.e. 0 is an
absorbing state, and it is verified $P(X_{0,1}=0)>0$ or $P(\phi_0(k)=0)>0$, $k=1,2,\ldots$,  it was established that  under A1) and A2),  if $\alpha  > \sigma^2/(2m)$ and an assumption on conditional moments holds, then $P(Z_n\to \infty) > 0.$ In the present paper we will consider critical CBPs, $\{Z_n\}_{n\geq 0}$, satisfying the above conditions,  but with a reflecting barrier at zero, namely, $P(\phi_n(0)>0)>0$. Thus $\{Z_n\}_{n\geq 0}$ will have a finite number of returns to the sate zero till the explosion to infinity, i.e.  $P(Z_n\to\infty)=1$.

Notice that under A1), $\varepsilon(k)=(k+\alpha)m^{-1}$, $k\geq 1$, and, for simplicity in the posterior calculations,  we will also assume throughout the paper that $\varepsilon(0)= \alpha m^{-1}$. }
{
\begin{nota} The controlled branching process we are considering is such that  migration may take place in the next generation no matter the size of the current generation (when there are no individuals in the populations only immigration is possible). BGW processes with immigration at 0 were considered firstly in \cite{foster} and \cite{pakes}.
\end{nota}
}

In next result we calculate the first and second moments of a CBP which verifies A1) and A2).

\begin{proposicion}\label{momentos} Let $\{Z_{n}\}_{n\geq 0}$ be a CBP
with  $E[Z_0^2]<\infty$ and satisfying hypotheses A1) and A2). It is verified as $k\to\infty$ that $$E[Z_k]=O(k) \mbox{ and } E[Z_k^2]=O(k^2).$$
\end{proposicion}

\noindent{\bf Proof.}
{From (\ref{esperanza}) and A1) it follows that

\begin{equation}\label{espe}
E[Z_{k}]=E[m\varepsilon(Z_{k-1})]=E[Z_0]+k\alpha,\quad k\geq 0.
\end{equation}

 Using (\ref{varianza}) we have
\begin{eqnarray*}
E[Var[Z_{k}\mid\mathcal{F}_{k-1}]]&=&  m^{-1}\sigma^2(E[Z_{0}]+(k-1)\alpha)+m^2
E[\nu^2(Z_{k-1})].
\end{eqnarray*}
Now,  from  A2), we have that there exists $C>0$ such that $\nu^2(k)\leq Ck$ for all $k>0$, so that $E[\nu^2(Z_{k-1})]\leq C E[Z_{n-1}]+\nu^2(0)$.
Consequently,  letting $M_1=3\max\{m^{-1}\sigma^2 E[Z_0],\ $ $Cm^2E[Z_0],\  m^2\nu^2(0)\}$ and $M_2=2\max\{m^{-1}\sigma^2\alpha, m^2C\alpha\}$, we have
\begin{equation}\label{varcond}
E[Var[Z_{k}\mid\mathcal{F}_{k-1}]]\leq  M_1+ M_2(k-1).
\end{equation}
Hence,
\begin{eqnarray*}
Var[Z_{k}]&=& E[Var[Z_{k}\mid\mathcal{F}_{k-1}]]+Var[E[Z_{k}\mid
\mathcal{F}_{k-1}]]\\
&\leq&  M_1+ M_2(k-1)+Var[Z_{k-1}]
\\& \leq & kM_1+2^{-1}M_2k(k-1)+Var[Z_0].
\end{eqnarray*}
The latter inequality proves that $E[Z_k^2]=O(k^2)$.}
\vspace*{0.25cm}

Next Lemma presents certain relationships among the random variables $\{X_{n,j}:\ n=0,1,\ldots;\ j=1,2,\ldots\}$ which can be easily verified.

\begin{lema} \label{lema1}Let $\left\{X_{n,j}: n=0,1,\ldots; j=1,2,\ldots\right\}$  be a sequence of
i.i.d. random variables
with mean $m=E[X_{0,1}]$ and $\sigma^2=Var[X_{0,1}]$, assumed finite.
Let denote $S_{k,l}=\sum_{j=1}^l (X_{k-1,j}-m)$, $k=1,2,\ldots$, $l=1,\ldots$, and for $j=1,\ldots, l$,  $ \tilde{S}_k^j(l)=\sum_{j^{\prime}\not=j}^{l}(X_{k-1,j^{\prime}}-m)$ and $M>0,\ M\in \mathbb{R}$.
It is verified that
$$E\left[\sum_{j=1}^l(X_{k-1,j}-m)^{2} \mathbb{I}_{\{|\tilde{S}_{k}^j(l)|>M\}}\right]\leq \frac{ l^2\sigma^4}{M^2}$$ and

$$E\left[\left(\sum_{j,j^{\prime},j\not=j^{\prime}}^l(X_{k-1,j}-m)(X_{k-1,j^{\prime}}-m)\right)^{2}\right]=2l(l-1)\sigma^4.$$
\end{lema}

\section{Main result}\label{sec:main}
We introduce for each
 $n\in\mathbb{N}$, a stochastic process $W_{n}(t) = n^{-1}Z_{\lfloor nt \rfloor}$, for
$t \geq 0$, $t\in \mathbb{R}$, $\lfloor \cdot \rfloor$ denoting the integer part. It is
easy to see that $\left\{W_{n}\right\}_{n\geq 1}$ is a sequence of random
functions that take values in $D_{[0,\infty)}[0,\infty)$, which is
the space of non-negative functions on $[0,\infty)$ that are right
continuous and have left limits. We also denote by
$C_c^{\infty}[0,\infty)$ the space of infinitely differentiable
functions on $[0,\infty)$ which have  a compact support. Throughout the paper ``$\stackrel{\mathcal{D}}{\to}$'' denotes the convergence
of random functions in the Skorokhod topology.
\begin{teorema}\label{teor:main}
Let $\{Z_{n}\}_{n\geq 0}$ be a CBP
with  $E[Z_0^2]<\infty$, satisfying hypotheses A1) and A2).
Then,
$W_n\stackrel{\mathcal{D}}{\rightarrow} W$,
as $n\to\infty,$ being $W$ a non-negative diffusion process, with
generator $T f(x) = \alpha f'(x) + \frac{1}{2}x
\sigma^{2}m^{-1}f''(x)$, for $f \in C^{\infty}_{c}[0,\infty)$.
The process $W$ is  the pathwise unique solution of the stochastic
differential equation
\begin{equation}\label{sde1}
\mathrm{d}W(t)=\alpha \mathrm{d}t + \sqrt{\sigma^2 m^{-1} (W(t)){^+}
}\mathrm{d}{\cal W}(t), \qquad t\geq 0,
\end{equation}
 with initial value $W(0)=0$,  denoting $x^+=\max\{x,0\}$, $x \in \mathbb{R}$, and where  ${\cal W}$ is a standard Wiener process.
\end{teorema}
\begin{nota}
Taking into account Theorem A\hspace*{-0.05cm}\ref{teo:exist} in Appendix, the stochastic differential equation (SDE) (\ref{sde1}) has a pathwise unique solution $\{X(t)^{(x)}\}_{t\geq 0}$ for all initial values $X(0)^{(x)}=x\in \mathbb{R}$. Moreover if $x\geq 0$, then $X(t)^{(x)}\geq 0$ almost surely for all $t\geq 0$.
\end{nota}

In order to prove  Theorem \ref{teor:main}, we will establish  previously the weak convergence of random step processes defined from a martingale difference created from the CBP.

We  introduce the following sequence of martingale differences $\{M_k\}_{k\geq 1}$ with respect the filtration $\{\mathcal{F}_k\}_{k\geq 0}$ as:
$$M_k= Z_k-E[Z_{k}\mid \mathcal{F}_{k-1}]= Z_k-Z_{k-1}-\alpha, \ k\geq 1.$$
Consider the random step processes:
\begin{equation}\label{eq:21}
\mathcal{M}_n(t)=\frac{1}{n}\left(Z_{0}+\sum_{k=1}^{\lfloor n t\rfloor} M_{k}\right)=\frac{1}{n} Z_{\lfloor nt \rfloor}-\frac{\lfloor nt \rfloor}{n}\alpha , \quad t \geq 0, \quad n \in \mathbb{N}.
\end{equation}
\begin{teorema}\label{teo:previo} Let $\{Z_{n}\}_{n\geq 0}$ be a CBP
with  $E[Z_0^2]<\infty$, satisfying hypotheses A1) and A2).
It is verified that
$$
\mathcal{M}_{n} \stackrel{\mathcal{D}}{\longrightarrow}\mathcal{M}, \quad \text { as } n \rightarrow \infty,
$$
where the limit process $\mathcal{M}$ is the pathwise unique solution of

\begin{equation}\label{sde2}
\mathrm{d}\mathcal{M}(t)=\sqrt{m^{-1}\sigma^{2}(\mathcal{M}(t)+\alpha t)^{+}}\mathrm{d} \mathcal{W}(t), \quad t \geq 0,\  \mbox{with initial value } \mathcal{M}(0)=0.
\end{equation}
\end{teorema}
{\bf Proof} As was done in \cite{bbp21}, we prove the result by applying  Theorem A\hspace*{-0.05cm}\ref{teor:ip} in Appendix with $\mathcal{U}=\mathcal{M},\  U_{n}(k)=n^{-1} M_{k},\  k \in \mathbb{N}$, $U_{n}(0)=n^{-1} Z_{0}$, $\mathcal{F}_{n}(k)=\mathcal{F}_{k}, k \geq 0$,
where $n \in \mathbb{N}$ (yielding $\mathcal{U}_{n}=\mathcal{M}_{n}, n \in \mathbb{N}$, as well), and with coefficient functions
$\beta: [0,\infty) \times \mathbb{R} \rightarrow \mathbb{R}$ and $\gamma: [0,\infty) \times \mathbb{R} \rightarrow \mathbb{R}$  given by
$$
\beta(t, x)=0, \quad \gamma(t, x)=\sqrt{m^{-1}\sigma^{2}\left(x+\alpha t\right)^{+}}, \quad t \geq 0, \quad x \in \mathbb{R}.
$$

 Firstly, we check that the SDE (\ref{sde2}) has a pathwise unique strong solution $\left\{\mathcal{M}(t)^{(x)}\right\}_{t\geq 0}$ for all initial values $\mathcal{M}(0)^{(x)}=x \in \mathbb{R}$. In fact, notice that if $\left\{\mathcal{M}(t)^{(x)}\right\}_{t\geq 0}$ is a strong solution of the SDE
(\ref{sde2}) with initial value $\mathcal{M}(0)^{(x)}=x \in \mathbb{R}$, then, by It\^o's formula, the process $\mathcal{P}(t)=\mathcal{M}(t)^{(x)}+\alpha t,$
$t \geq 0$, is a solution of the SDE
\begin{equation}\label{sde3}
\mathrm{d} \mathcal{P}(t)=\alpha \mathrm{d} t+\sqrt{m^{-1}\sigma^{2} \mathcal{P}(t)^{+}} \mathrm{d} \mathcal{W}(t), \quad t \geq 0, \  \mbox{with initial value } \mathcal{P}(0)=x.
\end{equation}
Conversely, if $\{\mathcal{P}(t)^{(x)}\}_{t\geq 0}$ is a strong solution of the SDE (\ref{sde3}) with initial value  $\mathcal{P}^{(x)}(0)=x\in \mathbb{R}$, then,
by It\^o's formula, the process $\mathcal{M}(t)=\mathcal{P}(t)^{(x)}-\alpha t,$
$t \geq 0$, is a strong solution of the SDE (\ref{sde2}) with initial value $\mathcal{M}(0)=x$. Notice that SDE (\ref{sde3}) is the same as SDE (\ref{sde1}), consequently the SDE (\ref{sde3}) and therefore the SDE (\ref{sde2}) as well admit a pathwise unique strong solution with arbitrary initial value, and
\begin{equation}\label{eq:igual}\{\mathcal{M}(t)+\alpha t\}_{t\geq 0}\stackrel{\mathcal{D}}{=}\{W(t)\}_{t\geq 0}.\end{equation}

 Let us see that $E\left[\left(U_{n}(k)\right)^{2}\right]<\infty$ for all $n=1,2\ldots$ and $k=0,1,2,\ldots$. Indeed, taking into account {(\ref{varcond})} in Proposition \ref{momentos},
{ \begin{equation}\label{eq:11}
  E\left[\left(U_{n}(k)\right)^{2}\right]=n^{-2} E\left[M_{k}^{2}\right]=E[Var[Z_k\mid\mathcal{F}_{k-1}]]\leq \frac{M_1+ M_2(k-1)}{n^{2}}<\infty
  \end{equation}}
   and, by the assumption in the statement of the theorem, $E\left[\left(U_{n}(0)\right)^{2}\right]=n^{-2} E\left[Z_{0}^{2}\right]<\infty,$ for $n=1,2,\ldots$
Moreover, $U_{n}(0)=n^{-1} Z_{0} \stackrel{\text { a.s. }}{\longrightarrow} 0$ as $n \rightarrow \infty$, especially $U_{n}(0) \stackrel{\mathcal{D}}{\longrightarrow} 0$ as $n \rightarrow \infty$.

For conditions (i), (ii) and (iii) of Theorem  A\hspace*{-0.05cm}\ref{teor:ip} in Appendix, we have to check that for each $T >0$, $T\in \mathbb{R}$, as $n\to\infty$:
\begin{itemize}
\item[a)] $\sup _{t \in[0, T]}\left|\frac{1}{n} \sum_{k=1}^{\lfloor n t\rfloor} E\left[M_{k} \mid \mathcal{F}_{k-1}\right]-0\right| \stackrel{\mathrm{P}}{\longrightarrow} 0.$
\item[b)]$\sup _{t \in[0, T]} \left|\frac{1}{n^{2}} \sum_{k=1}^{\lfloor n t\rfloor} E\left[M_{k}^{2} \mid \mathcal{F}_{k-1}\right]-\int_{0}^{t} \frac{\sigma^{2}}{m}\left(\mathcal{M}_{n}(s)+\alpha s\right)^{+} \mathrm{d} s\right| \stackrel{\mathrm{P}}{\longrightarrow} 0.$
\item[c)] For all $\theta>0, \theta \in \mathbb{R}$,
$\frac{1}{n^{2}} \sum_{k=1}^{\lfloor n T\rfloor} E\left[M_{k}^{2} \mathbb{I}_{\left\{\left|M_{k}\right|>n \theta\right\}} \mid \mathcal{F}_{k-1}\right] \stackrel{\mathrm{P}}{\longrightarrow} 0.$
 \end{itemize}
 Since $E\left[M_{k} \mid \mathcal{F}_{k-1}\right]=0$, $n \in \mathbb{N}$, $k \in \mathbb{N}$, a) holds.

Let us check b).

\noindent For each $s>0, s \in \mathbb{R}$ and $n \in \mathbb{N}$, and  for all $t >0$, $t\in \mathbb{R}$ and $n\in\mathbb{N}$, we have:
$$
\mathcal{M}_{n}(s)+\alpha s=\frac{1}{n} Z_{\lfloor n s\rfloor}+\frac{n s-\lfloor n s\rfloor }{n} \alpha,
$$
thus $\left(\mathcal{M}_{n}(s)+\alpha s\right)^{+}=\mathcal{M}_{n}(s)+\alpha s$. Now, we have, for all $t>0$ and $n \in \mathbb{N}$,
\begin{eqnarray*}
\int_{0}^{t}\left(\mathcal{M}_{n}{(s)}+\alpha s\right)^{+} \mathrm{d} s&=&\int_{0}^{t}\left(\frac{1}{n} Z_{\lfloor n s\rfloor}+\frac{n s-\lfloor n s\rfloor}{n} \alpha\right) \mathrm{d} s
\\&=&\sum_{k=0}^{\lfloor n t\rfloor -1} \int_{k / n}^{(k+1) / n}\left(\frac{1}{n} Z_{k}+\frac{n s-k}{n} \alpha\right) \mathrm{d} s\\&&+\int_{\lfloor n t\rfloor / n}^{t}\left(\frac{1}{n} Z_{\lfloor n t \rfloor}+\frac{n s-\lfloor n t\rfloor}{n} \alpha\right) \mathrm{d} s\\
&=&\frac{1}{n^{2}} \sum_{k=0}^{\lfloor n t \rfloor-1} Z_{k}+\frac{n t-\lfloor n t\rfloor}{n^{2}} Z_{\lfloor n t\rfloor}+\frac{\alpha}{2n^2}\lfloor n t \rfloor\\&&+\frac{\alpha}{n}\left(\frac{n}{2}\left(t^{2}-\frac{\lfloor n t\rfloor^{2}}{n^{2}}\right)-\lfloor n t\rfloor\left(t-\frac{\lfloor n t\rfloor}{n}\right)\right)\\
&=&\frac{1}{n^{2}} \sum_{k=0}^{\lfloor n t \rfloor-1} Z_{k}+\frac{n t-\lfloor n t\rfloor}{n^{2}} Z_{\lfloor n t \rfloor}+\frac{\lfloor n t\rfloor+(n t-\lfloor n t\rfloor)^{2}}{2 n^{2}} \alpha. \end{eqnarray*}

It is verified that, for $t>0$ and $n\in\mathbb{N}$,
\begin{eqnarray*}
\frac{1}{n^{2}} \sum_{k=1}^{\lfloor n t\rfloor} E\left[M_{k}^{2} \mid \mathcal{F}_{k-1}\right]&=&\frac{1}{n^{2}} \sum_{k=1}^{\lfloor n t\rfloor}Var[Z_k\mid \mathcal{F}_{k-1}]\\
&=&\frac{1}{n^{2}} \sum_{k=1}^{\lfloor n t\rfloor} \left(m^2\nu^2(Z_{k-1})+\frac{\sigma^2}{m}(Z_{k-1}+\alpha)\right)\\
& =&\frac{m^2}{n^2} \sum_{k=1}^{\lfloor n t\rfloor}\nu^2(Z_{k-1})+ \frac{\lfloor n t\rfloor \alpha \sigma^2}{n^2m}+ \frac{\sigma^2}{n^2 m} \sum_{k=1}^{\lfloor n t\rfloor}Z_{k-1}.
\end{eqnarray*}

Consequently,
\begin{eqnarray*}
\frac{1}{n^{2}} \sum_{k=1}^{\lfloor n t\rfloor} E\left[M_{k}^{2} \mid \mathcal{F}_{k-1}\right]&-&\int_{0}^{t} \frac{\sigma^{2}}{m}\left(\mathcal{M}_{n}(s)+\alpha s\right)^{+} \mathrm{d} s= \frac{m^2}{n^2} \sum_{k=1}^{\lfloor n t\rfloor}\nu^2(Z_{k-1})+\frac{\lfloor n t\rfloor \alpha \sigma^2}{n^2m}\\
&&- \frac{\sigma^2 (nt -\lfloor nt \rfloor)}{m n^2} Z_{\lfloor nt \rfloor}-\frac{\sigma^2}{m}\frac{\lfloor n t\rfloor+(n t-\lfloor n t\rfloor)^{2}}{2 n^{2}} \alpha.
\end{eqnarray*}
Since for each $T>0, T \in \mathbb{R}$,
$$
\begin{array}{l}
\sup _{t \in[0, T]} \frac{\lfloor n t\rfloor}{n^{2}} \leqslant \frac{T}{n} \rightarrow 0 ,\quad \text { as } n \rightarrow \infty \\[0.1in]
\sup _{t \in[0, T]} \frac{\lfloor n t\rfloor+(n t-\lfloor n t\rfloor)^{2}}{2 n^{2}} \leqslant \frac{T}{2 n}+\frac{1}{2 n^{2}} \rightarrow 0, \quad \text { as } n \rightarrow \infty,
\end{array}
$$
in order to show $b)$, it suffices to prove that for each $T>0, T \in \mathbb{R}$,
\begin{equation}\label{eq:28}
\frac{1}{n^{2}} \sup _{t \in[0, T]}\left((n t-\lfloor n t\rfloor) Z_{\lfloor n t\rfloor}\right) \leq \frac{1}{n^{2}} \sup _{t \in[0, T]} Z_{\lfloor n t\rfloor} \stackrel{{P}}{\longrightarrow} 0 \quad \mbox{ as } n \rightarrow \infty .
\end{equation}
  and
 \begin{equation}\label{eq:28b}
\frac{m^2}{n^{2}} \sup _{t \in[0, T]} \sum_{k=1}^{\lfloor n t\rfloor}\nu^2(Z_{k-1}) \stackrel{{P}}{\longrightarrow} 0 \quad \mbox{ as } n \rightarrow \infty .
\end{equation}

First we check (\ref{eq:28}). For each $k \in \mathbb{N}$, we have $Z_{k}=Z_{k-1}+M_{k}+\alpha$, thus
$$
Z_{k}=Z_{0}+\sum_{j=1}^{k} M_{j}+k \alpha,$$
and hence, for each $t >0,\ t \in \mathbb{R}$ and $n \in \mathbb{N}$, we get
$$
Z_{\lfloor n t]}=\left|Z_{\lfloor n t\rfloor}\right| \leqslant Z_{0}+\sum_{j=1}^{\lfloor n t\rfloor}\left|M_{j}\right|+\lfloor n t\rfloor \alpha.
$$
Consequently, in order to prove $(\ref{eq:28})$, it suffices to show
$$
\frac{1}{n^{2}} \sup _{t \in[0, T]} \sum_{j=1}^{\lfloor n t\rfloor}\left|M_{j}\right| \leqslant \frac{1}{n^{2}} \sum_{j=1}^{\lfloor n T\rfloor}\left|M_{j}\right| \stackrel{{P}}{\longrightarrow} 0, \quad \text { as } n \rightarrow \infty.
$$
By (\ref{eq:11}), $E[M_k^2]=O(k),$ as $k\to\infty$, and therefore by Jensen's inequality, $E[|M_k|]=O(k^{1/2}),$ as $k\to\infty$, and hence
$$
E\left[\frac{1}{n^{2}} \sum_{j=1}^{\lfloor n T\rfloor}\left|M_{j}\right|\right]=\frac{1}{n^{2}} \sum_{j=1}^{\lfloor n T\rfloor} {O}\left(j^{1 / 2}\right)={O}\left(n^{-1 / 2}\right) \rightarrow 0, \quad \text { as } n \rightarrow \infty.
$$
Thus we obtain $n^{-2} \sum_{j=1}^{\lfloor n T\rfloor}\left|M_{j}\right| \stackrel{{P}}{\longrightarrow} 0$ as $n \rightarrow \infty$ implying $(\ref{eq:28})$.

Now, taking into account  hypothesis A2)
\begin{equation}\label{eq:12}E\left[\frac{1}{n^2}\sum_{k=1}^{\lfloor n T\rfloor}\nu^2(Z_{k-1})\right]=\frac{1}{n^2}\sum_{k=1}^{\lfloor n T\rfloor}O(k^{\beta})=O(n^{\beta-1}), \end{equation}
and hence (\ref{eq:28b}).

Let us check c).

We write
$$M_k=\sum_{j=1}^{\phi_{k-1}(Z_{k-1})}(X_{k-1,j}-m)+m(\phi_{k-1}(Z_{k-1})-\varepsilon(Z_{k-1})).$$
Let denote $N_k=\sum_{j=1}^{\phi_{k-1}(Z_{k-1})}(X_{k-1,j}-m)$. It is verified for each $n,k\in \mathbb{N}$, and $\theta>0, \theta\in \mathbb{R}$ that
$$M_k^2\leq 2\left(N_k^2+ m^2(\phi_{k-1}(Z_{k-1})-\varepsilon(Z_{k-1}))^2\right),$$  and
$$\mathbb{I}_{\left\{\left|M_{k}\right|>n \theta\right\}}\leq \mathbb{I}_{\left\{\left|N_{k}\right|>n \theta/2\right\}}+\mathbb{I}_{\left\{\left|\phi_{k-1}(Z_{k-1})-\varepsilon(Z_{k-1})\right|>n \theta/2m\right\}}.$$ Hence
$$M_k^2\mathbb{I}_{\left\{\left|M_{k}\right|>n \theta\right\}}\leq 2N_k^2\mathbb{I}_{\left\{\left|N_{k}\right|>n \theta/2\right\}}+2N_k^2\mathbb{I}_{\left\{\left|\phi_{k-1}(Z_{k-1})-\varepsilon(Z_{k-1})\right|>n \theta/2m\right\}}+2m^2(\phi_{k-1}(Z_{k-1})-\varepsilon(Z_{k-1}))^2.
$$

In consequence, to check c) we will prove, as $n\to\infty$,
\begin{itemize}
\item[c.1)]
$\frac{1}{n^{2}} \sum_{k=1}^{\lfloor nT\rfloor} E\left[N_{k}^{2} \mathbb{I}_{\left\{\left|N_{k}\right|>n \theta\right\}} \mid \mathcal{F}_{k-1}\right] \stackrel{{P}}{\longrightarrow} 0 \quad$ for all $\theta>0, \theta \in \mathbb{R}$.
\item[c.2)]
$\frac{1}{n^{2}} \sum_{k=1}^{\lfloor nT\rfloor} E\left[N_{k}^{2} \mathbb{I}_{\left\{\left|\phi_{k-1}(Z_{k-1})-\varepsilon(Z_{k-1})\right|>n \theta\right\}} \mid \mathcal{F}_{k-1}\right] \stackrel{{P}}{\longrightarrow}0$ for all $\theta>0, \theta \in \mathbb{R}$.
\item[c.3)]
$\frac{1}{n^{2}} \sum_{k=1}^{\lfloor nT\rfloor} E\left[(\phi_{k-1}(Z_{k-1})-\varepsilon(Z_{k-1}))^{2} \mid \mathcal{F}_{k-1}\right]\stackrel{{P}}{\longrightarrow} 0.
$
\end{itemize}
In what follows let $\theta>0, \theta \in \mathbb{R}$ be fixed.

Let us see c.3). It is verified that
\begin{eqnarray*}
\frac{1}{n^{2}} \sum_{k=1}^{\lfloor nT\rfloor} E\left[(\phi_{k-1}(Z_{k-1})-\varepsilon(Z_{k-1}))^{2} \mid \mathcal{F}_{k-1}\right]&= &\frac{1}{n^{2}} \sum_{k=1}^{\lfloor nT\rfloor}Var\left[\phi_{k-1}(Z_{k-1})\mid \mathcal {F}_{k-1}\right]\\ &=&
\frac{1}{n^{2}} \sum_{k=1}^{\lfloor nT\rfloor} \nu^2(Z_{k-1})\stackrel{{P}}{\longrightarrow} 0, \mbox{ as } n\to\infty.
\end{eqnarray*}
This latter was proved by considering (\ref{eq:12}).

Now, we check c.1).
By the properties of conditional expectation with respect to a $\sigma$ -algebra, we get for all $n, k \in \mathbb{N}$,
$$
E\left[N_{k}^{2} \mathbb{I}_{\left\{\left|N_{k}\right|>n \theta\right\}} \mid \mathcal{F}_{k-1}\right]=F_{n, k}\left(Z_{k-1}\right),
$$
where on $\{Z_{k-1}=z\}$, with $z=0,1,\ldots$
$$
F_{n, k}(z)=E\left[S_{k}(z)^{2} \mathbb{I}_{\left\{\left|S_{k}(z)\right|>n \theta\right\}}\right], \mbox{ where }
S_{k}(z)=\sum_{j=1}^{\phi_{k-1}(z)}\left(X_{k-1,j}-m\right).$$
Consider the decomposition $F_{n, k}(z)=A_{n, k}(z)+B_{n, k}(z)$ with
\label{page}
\begin{eqnarray*}
A_{n, k}(z)&=&E\left[\sum_{j=1}^{\phi_{k-1}(z)} (X_{k-1,j}-m)^{2} \mathbb{I}_{\{|S_{k}(z)|>n \theta\}}\right], \\
B_{n, k}(z)&=&E\left[\sum_{j, j^{\prime},j\not=j^{\prime}}^{\phi_{k-1}(z)} (X_{k-1,j}-m)(X_{k-1,j^{\prime}}-m) \mathbb{I}_{\{\left|S_{k}(z)\right|>n \theta\}}\right].
\end{eqnarray*}

Now, let denote $S_{k,l}=\sum_{j=1}^l (X_{k-1,j}-m)$, $k=1,2,\ldots$, $l=0,1,\ldots$. It is verified the inequality, for $j\in \{1,\ldots, l\}$
$$|S_{k,l}|\leq |X_{k-1,j}-m|+|\tilde{S}_k^j(l)|,\mbox{ with } \tilde{S}_k^j(l)=\sum_{j^{\prime}\not=j}^{l}(X_{k-1,j^{\prime}}-m).$$
We have, using Lemma \ref{lema1},
\begin{eqnarray*}
A_{n, k}(z)&=&E\left[E\left[\left.\sum_{j=1}^{\phi_{k-1}(z)} (X_{k-1,j}-m)^{2} \mathbb{I}_{\{|S_{k}(z)|>n \theta\}}\right|\phi_{k-1}(z)\right]\right]\\&=&
\sum_{l=0}^{\infty}E\left[\sum_{j=1}^{l} (X_{k-1,j}-m)^{2} \mathbb{I}_{\{|{S}_{k,l}|>n \theta\}}\right]P(\phi_{k-1}(z)=l)\\&\leq&
\sum_{l=0}^{\infty}\sum_{j=1}^{l} (E\left[(X_{k-1,j}-m)^{2} \mathbb{I}_{\{|X_{k-1,j}-m|>n \theta/2\}}\right]\\&&+ \left.E[(X_{k-1,j}-m)^{2} \mathbb{I}_{\{|\tilde{S}_{k}^j(z)|>n \theta/2\}}]\right)P(\phi_{k-1}(z)=l)\\
&\leq & \sum_{l=0}^{\infty} \left(lE\left[(X_{0,1}-m)^{2} \mathbb{I}_{\{|X_{0,1}-m|>n \theta/2\}}\right]+ \frac{4l^2\sigma^4}{n^2\theta^2}\right)P(\phi_{k-1}(z)=l)\\
&=&\varepsilon(z)E\left[(X_{0,1}-m)^{2} \mathbb{I}_{\{|X_{0,1}-m|>n \theta/2\}}\right]+\frac{4\sigma^4}{n^2\theta^2}E[(\phi_{0}(z))^2]
\end{eqnarray*}
Therefore
$$A_{n,k}(z)\leq A_{n,k}(z)^{(1)} +A_{n,k}(z)^{(2)},$$
with
\begin{eqnarray*}
A_{n,k}(z)^{(1)}&=& \varepsilon(z)E\left[(X_{0,1}-m)^{2} \mathbb{I}_{\{|X_{0,1}-m|>n \theta/2\}}\right],\\
A_{n,k}(z)^{(2)}&=& (\nu^2(z)+(\varepsilon(z))^2)\frac{4\sigma^4}{n^2\theta^2}.
\end{eqnarray*}

{Using (\ref{espe}) in Proposition \ref{momentos}}, it is verified that for $n\in\mathbb{N}$
\begin{eqnarray*}
E\left[\frac{1}{n^2}\sum_{k=1}^{\lfloor nT\rfloor} A_{n,k}(Z_{k-1})^{(1)}\right]&=&
\frac{1}{n^2}\sum_{k=1}^{\lfloor nT\rfloor}E[\varepsilon(Z_k)]E\left[(X_{0,1}-m)^{2} \mathbb{I}_{\{|X_{0,1}-m|>n \theta/2\}}\right]\\
&=&\frac{1}{n^2}\left(\sum_{k=1}^{\lfloor nT\rfloor}O(k)\right)E\left[(X_{0,1}-m)^{2} \mathbb{I}_{\{|X_{0,1}-m|>n \theta/2\}}\right]\\&=&
O(1)E\left[(X_{0,1}-m)^{2} \mathbb{I}_{\{|X_{0,1}-m|>n \theta/2\}}\right].
\end{eqnarray*}
By applying the dominated convergence theorem we have
\begin{equation}\label{eq:conv1}
E\left[\frac{1}{n^2}\sum_{k=1}^{\lfloor nT\rfloor} A_{n,k}(Z_{k-1})^{(1)}\right]\to 0,\ \mbox{as }n\to\infty.
\end{equation}
It is also verified by using { again (\ref{espe}) in Proposition \ref{momentos} and A2)} that
\begin{eqnarray}
E\left[\frac{1}{n^2}\sum_{k=1}^{\lfloor nT\rfloor} A_{n,k}(Z_{k-1})^{(2)}\right]&=&\frac{1}{n^2}\sum_{k=1}^{\lfloor nT\rfloor}E\left[\left(\nu^2(Z_{k-1})+(\varepsilon(Z_{k-1}))^2\right)\frac{4\sigma^4}{n^2\theta^2}\right]\nonumber
\\&=&\frac{4\sigma^4}{n^4\theta^2}\sum_{k=1}^{\lfloor nT\rfloor}O(k^2)=O(n^{-1})\to 0, \mbox{ as } n\to\infty.\label{eq:conv2}
\end{eqnarray}
Taking into account (\ref{eq:conv1}) and (\ref{eq:conv2})
we have that, as $n\to\infty$,
$$\frac{1}{n^2}\sum_{k=1}^{\lfloor nT\rfloor} A_{n,k}(Z_{k-1})\stackrel{{P}}{\longrightarrow} 0.$$
Let us now dealt with $B_{n,k}(z)$. It is verified that, using Cauchy-Schwarz's inequality:
\begin{eqnarray*}
B_{n,k}(z)&=&\sum_{l=0}^{\infty}E\left[\sum_{j,j^{\prime},j\not=j^{\prime}}^l(X_{k-1,j}-m)(X_{k-1,j^{\prime}}-m)\mathbb{I}_{\{|{S}_{k,l}|>n \theta\}}\right]P(\phi_{k-1}(z)=l)\\
&\leq& \sum_{l=0}^{\infty}\sqrt{E\left[\left(\sum_{j,j^{\prime},j\not=j^{\prime}}^l(X_{k-1,j}-m)(X_{k-1,j^{\prime}}-m)\right)^2\right]E[\mathbb{I}_{\{|{S}_{k,l}|>n \theta\}}]}P(\phi_{k-1}(z)=l).
\end{eqnarray*}
Now, using Markov's inequality
$$E[\mathbb{I}_{\{|{S}_{k,l}|>n \theta\}}]\leq \frac{Var[S_{k,l}]}{n^2\theta^2}=\frac{l\sigma^2}{n^2\theta^2},$$ and using
Lemma \ref{lema1}, we have
$$B_{n,k}(z)\leq \sum_{l=0}^\infty \sqrt{2l^2\sigma^4n^{-2}\theta^{-2}l\sigma^2}P(\phi_{k-1}(z)=l)=\frac{\sqrt{2}\sigma^3}{\theta n}E[(\phi_{0}(z))^{3/2}].$$

By Lyapunov's inequality, $E[(\phi_{k-1}(z))^{3/2}]^{2/3}\leq E[(\phi_{k-1}(z))^{2}]^{1/2}=(\nu^2(z)+(\varepsilon(z))^2)^{1/2}$. Consequently, in order to prove, as $n\to\infty$,
$$\frac{1}{n^2}\sum_{k=1}^{\lfloor nT\rfloor} B_{n,k}(Z_{k-1})\stackrel{{P}}{\longrightarrow} 0,$$
is enough to check
that
$$\frac{1}{n^3}\sum_{k=1}^{\lfloor nT\rfloor} (\nu^2(Z_{k-1})+(\varepsilon(Z_{k-1}))^2)^{3/4}\stackrel{{P}}{\longrightarrow} 0.$$
In fact, using  hypotheses A1) and A2) and Proposition \ref{momentos}, we have
$$E\left[\frac{1}{n^3}\sum_{k=1}^{\lfloor nT\rfloor} (\nu^2(Z_{k-1})+(\varepsilon(Z_{k-1}))^2)^{3/4}\right]=n^{-3}\sum_{k=1}^{\lfloor nT\rfloor}O(k^{3/2})=O(n^{-1/2}).$$

Finally, we check c.2). We have that
$$E\left[N_{k}^{2} \mathbb{I}_{\left\{\left|\phi_{k-1}(Z_{k-1})-\varepsilon(Z_{k-1})\right|>n \theta\right\}} \mid \mathcal{F}_{k-1}\right]=G_{n,k}(Z_{k-1}),$$
where on $\{Z_{k-1}=z\}$, with $z=0,1,\ldots$,
$$G_{n,k}(z)=E[S_k(z)^2\mathbb{I}_{\{|\phi_{k-1}(z)-\varepsilon(z)|>n\theta\}}].$$
Now, again by Cauchy-Schwarz's inequality and Markov's inequality
\begin{eqnarray*}
G_{n,k}(z)&=&\sum_{l=0}^\infty \mathbb{I}_{\{|l-\varepsilon(z)|>n\theta\}}E[S_{k,l}^2]P(\phi_{k-1}(z)=l)=\sigma^2E[\phi_{k-1}(z)\mathbb{I}_{\{|\phi_{k-1}(z)-\varepsilon(z)|>n\theta\}}]\\
 &\leq&\sigma^2\sqrt{E[\phi_{k-1}^2(z)]P(|\phi_{k-1}(z)-\varepsilon(z)|>n\theta)}\leq \sigma^2E[(\phi_{0}(z))^2]^{1/2}\left(\frac{\nu^2(z)}{n^2\theta^2}\right)^{1/2}.
\end{eqnarray*}
In consequence from
$$E\left[\frac{\sigma^2}{\theta n^3}\sum_{k=1}^{\lfloor nT\rfloor}\left(\nu^2(Z_{k-1})+(\varepsilon(Z_{k-1}))^2\right)^{1/2}(\nu^2(Z_{k-1}))^{1/2}\right]=O(n^{\beta/2-1}),$$
c.2) follows.

\vspace*{1cm}

Finally, using the weak convergence of $\{\mathcal{M}_n\}_{n\geq 1}$, we will obtain weak convergence of $\{W_n\}_{n\geq 1}$.
\vspace*{.25cm}

\noindent{\bf Proof of Theorem \ref{teor:main}.}
A version of the continuous mapping theorem is applied (see Lemma \ref{lema:cmt} in Appendix). {Let $D_{\mathbb{R}}[0,\infty)$ be the space of the real functions on $[0,\infty)$ that are right
continuous and have left limits}. For each $n \in \mathbb{N}$, by (\ref{eq:21}), $\{W_n(t)\}_{t \geq 0}=\Psi^{(n)}\left(\mathcal{M}_{n}\right)$, where the mapping
{$\Psi^{(n)}: D_{\mathbb{R}}[0,\infty)  \rightarrow D_{\mathbb{R}}[0,\infty) $} is given by
$$
\left(\Psi^{(n)}(f)\right)(t)=f\left(\frac{\lfloor n t\rfloor}{n}\right)+\frac{\lfloor n t\rfloor}{n} \alpha,
$$
for $f \in D_{\mathbb{R}}[0,\infty) $ and $t \in [0,\infty)$. Indeed, for each $n\in\mathbb{N}$ and $t\geq 0$:
$$ \left(\Psi^{(n)}(\mathcal{M}_n)\right)(t)=\mathcal{M}_n(\lfloor nt \rfloor/n)+\frac{\lfloor nt \rfloor}{n}\alpha=\frac{1}{n}Z_{\lfloor nt \rfloor}=W_n(t).$$

Further, taking into account (\ref{eq:igual}), $W \stackrel{\mathcal{D}}{=} \Psi(\mathcal{M})$, where the mapping $\Psi: D_{\mathbb{R}}[0,\infty) \rightarrow D_{\mathbb{R}}[0,\infty)$ is given by
$$
(\Psi(f))(t)=f(t)+\alpha t, \quad f \in D_{\mathbb{R}}[0,\infty), \quad t \in [0,\infty).
$$
The measurability of $\Psi^{(n)}$, $n\in\mathbb{N}$ and $\Psi$  likewise the conditions for applying Lemma \ref{lema:cmt} are proved in \cite{bbp21}.
\begin{nota}

\noindent 1) Notice that the result in Theorem \ref{teor:main} is also valid as $\alpha=0$, and even when the hypothesis A1) is replaced with the more general condition $\tau_m(k) = 1 + k^{-1}\alpha+o(k^{-1})$, as $\ k \to\infty$, $\alpha\geq 0$, being the calculation in this latter scenario a little more cumbersome.  In the  case $\alpha=0$ and $m=1$ the result  provides an alternative proof of the weak convergence result for the BGW process (see \cite{ek}, p. 388) for {a non-array version}.
\vspace*{0.2cm}

\noindent {2) As noted in the introduction, a BPI can be written as a special case of a CBP. For this case,  by considering  $m = E[X_{0,1}]=1$ ($\varepsilon(k)=k+E[I_0]$ and $\nu^2(k)=Var[I_0]$), we obtain an analogous result to that in \cite{bbp21}. Recall that the difference between both models is in which generation the immigrants will give rise to their offspring.}

\vspace*{0.2cm}

\noindent {3) As was pointed out in the Introduction, the proof of the main result follows similar steps as those  given in \cite{bbp21}. One can check that  similar formulas often appear being the roles of the immigration  mean and the offspring variance in the BPI case   played in the CBP by $\alpha$ and $m^{-1}\sigma^2$, respectively. However, new approaches  by considering conditioning arguments are needed to dealt with b)- c) in Theorem \ref{teor:main}, as a consequence that random sums of i.i.d. random variables arise in the proofs, see for instance the definition of $A_{n,k}(z)$ and $B_{n,k}(z)$, in p.\pageref{page}. An extra work is required to calculate the mathematical expectation of these quantities.}
\end{nota}
\noindent{\bf Acknowledgements:} The authors thank Professor M. Barczy for his constructive suggestions which have improved this paper.
This research has been supported by the Ministerio de Ciencia e Innovaci\'on of Spain (grant PID2019-108211GB-I00/AEI/10.13039/501100011033).
\section*{Appendix}

\begin{teoremaA}\label{teor:ip}
Let $\beta: [0,\infty) \times \mathbb{R} \rightarrow \mathbb{R}$ and $\gamma: [0,\infty) \times \mathbb{R} \rightarrow \mathbb{R}$ be continuous functions. Assume that uniqueness in the sense of probability law holds for the SDE
\begin{equation}\label{C.1}
\mathrm{d} \mathcal{U}(t)=\beta\left(t, \mathcal{U}(t)\right) \mathrm{d} t+\gamma\left(t, \mathcal{U}(t)\right) \mathrm{d} \mathcal{W}(t), \quad t \geq 0,
\end{equation}
with initial value $\mathcal{U}(0)=u(0)$ for all $u(0) \in \mathbb{R}$, where $\mathcal{W}=\left\{\mathcal{W}(t)\right\}_{t\geq 0}$
is an one-dimensional standard Wiener process. Let $\mathcal{U}=\left\{\mathcal{U}(t)\right\}_{t\geq 0}$ be a solution of $(\ref{C.1})$ with initial value $\mathcal{U}(0)=0.$
For each $n \in \mathbb{N}$, let $\left\{U_{n}(k): k=0,1,2,\ldots\right\}$ {be a sequence of real-valued random variables adapted to a filtration $\left\{\mathcal{F}_{n}(k):  k=0,1,2,\ldots\right\}$, that is, $U_{n}(k)$ is $ \mathcal{F}_{n}(k)$- measurable}. Let
$$
\mathcal{U}_{n}(t):=\sum_{k=0}^{\lfloor n t\rfloor} U_{n}(k), \quad t\geq 0, \quad n \in \mathbb{N}.
$$
Suppose that $E\left[\left(U_{n}(k)\right)^{2}\right]<\infty$ for all $n, k \in \mathbb{N}$, and $\mathcal{U}_n(0) \stackrel{\mathcal{D}}{\longrightarrow} 0$ as $n \rightarrow \infty$. Suppose that
for each $T >0$
\begin{itemize}
\item[(i)] $\sup _{t \in[0, T]}\left|\sum_{k=1}^{\lfloor n t\rfloor} E\left[U_{n}(k) \mid \mathcal{F}_{n}(k-1)\right]-\int_{0}^{t} \beta\left(s, \mathcal{U}_{n}(s)\right) \mathrm{d} s\right| \stackrel{{P}}{\longrightarrow} 0$ as $n \rightarrow \infty$,
\item[(ii)] $\sup _{t \in[0, T]}\left|\sum_{k=1}^{\lfloor n t\rfloor} \operatorname{Var}\left[U_{n}(k) \mid \mathcal{F}_{n}(k-1)\right]-\int_{0}^{t}\left(\gamma\left(s, \mathcal{U}_{n}(s)\right)\right)^{2} \mathrm{~d}s\right| \stackrel{{P}}{\longrightarrow} 0$ as $n \rightarrow \infty$,
\item[(iii)] $\sum_{k=1}^{\lfloor n T\rfloor} E\left[\left(U_{n}(k)\right)^{2} {\mathbb{I}}_{\left\{\left|U_{n}(k)\right|>\theta\right\}} \mid \mathcal{F}_{n}(k-1)\right] \stackrel{{P}}{\longrightarrow} 0$ as $n \rightarrow \infty$ for all $\theta >0$.
\end{itemize}
Then $\mathcal{U}_{n} \stackrel{\mathcal{D}}{\longrightarrow} \mathcal{U}$ as $n \rightarrow \infty$.
\end{teoremaA}
The proof can be seen in \cite{ip10}.

\begin{teoremaA}\label{teo:exist} Let $a$, $b$, $c$ real constants such that $a>0$. Consider the stochastic differential equation
\begin{equation}\label{eq:sde1}
\mathrm{d}X(t)=(bX(t)+c)\mathrm{d}t +\sqrt{2aX(t)^+}\mathrm{d}\mathcal{W}_t, \ t\geq 0.
\end{equation}
There exists a pathwise unique strong solution $\{X(t)^{(x)}\}_{t\geq 0}$ for all initial values $X(0)^{(x)}=x\in \mathbb{R}$. Moreover if $x\geq 0$, then $X(t)^{(x)}\geq 0$ almost surely for all $t\geq 0$.
In the case $c\geq 0$, the solution of (\ref{eq:sde1}) defines  diffusion process with generator
$$Tf(x)=(bx+c)f'(x) + axf''(x),\ f\in C_c^\infty[0,\infty),$$
where $C_c^\infty[0,\infty)$ is  the space of infinitely differentiable functions on $[0,\infty)$ which have a compact
support.
\end{teoremaA}
The proof can be seen in \cite{iw89} p. 235.

\begin{lemaA}\label{lema:cmt}
Let $S$ and $T$ be two metric spaces, and $X, X_1, X_2,
\cdots$ be random functions with values in $S$ with $X_n
\stackrel{\mathcal{D}}{\rightarrow} X$. Consider some measurable mappings $h,
h_1, h_2, \cdots: S \rightarrow T$ and a measurable set $C \subset S$
with $X \in C$ a.s. such that $h_n(s_n) \rightarrow h(s)$ as $s_n
\rightarrow s \in C$. Then $h_n(X_n) \stackrel{d}{\rightarrow}
h(X)$.
\end{lemaA}
The previous version of the continuous mapping theorem can be found in Theorem 3.27 in \cite{kale}.


\begin{thebibliography}{99}

\bibitem{bbp21} Barczy, M., Bezdány, D. and Pap, G. A note on asymptotic behavior of critical Galton-Watson processes with immigration, arXiv preprint, arXiv:2103.07878v2:1--20, 2021.

\bibitem{ek} Ethier, S.N., Kurtz, T.G. {\it Markov Processes: Characterization and
Convergence.} Wiley Series in Probability and Mathematical Statistics. John Wiley \& Sons, 1986.

\bibitem{feller} Feller, W.  Diffusion processes in genetics. In: Proc Second Berkeley
Symp Math Statist Prob, University of California Press, Berkeley, 227-
246, 1951.

\bibitem{foster} {Foster, J.H.  A limit theorem for a branching process with state-dependent immigration. {\it The Annals of Mathematical Statistics},
42: 1773-1776, 1971.}

 \bibitem{b2}
Gonz\'alez, M., Molina, M., and del Puerto, I. Asymptotic behaviour for the
  critical controlled branching process with random control function.
{\it Journal of Applied Probability}, 42: 463--477, 2005.

 \bibitem{b3} González, M. and del Puerto, I. Diffusion approximation of an array of controlled branching processes. {\it Methodology and Computing in Applied Probability,} 14: 843-861, 2012.

\bibitem{gpy}  González, M., del Puerto, I.  and Yanev, G.P. {\it Controlled Branching Processes}. ISTE Ltd and
John Wiley and Sons, Inc., 2018.

\bibitem{iw89} Ikeda, N. and Watanabe, S. . Stochastic {\it Differential Equations and Diffusion
Processes}, 2nd ed. North-Holland, Kodansha, Amsterdam, Tokyo, 1989.

\bibitem{ip10} Ispány, M. and Pap, G.  A note on weak convergence of step processes. {\it Acta Mathematica Hungarica}, 126(4): 381-395, 2010.


\bibitem{jirina} Jirina, M.  Stochastic Branching Processes with Continuous
State Space. {\it Czechoslovak Mathematical Journal}, 8:  292--313, 1958.

\bibitem{kale}
{Kallenberg, O.} {\it Foundations of Modern Probability}. Springer, 1997

\bibitem{lindvall} Lindvall, T. Convergence of critical Galton-Watson processes. {\it Journal of Applied Probability}, 9: 445-450, 1972.

\bibitem{pakes} {Pakes, A.G. On the critical Galton-Watson process with immigration. {\it
 Journal of the Australian Mathematical Society}, 12: 476-482, 1971.}

\bibitem{sriram} Sriram, T.N.,   Bhattacharya, A.,  González, M., Martínez, R., del Puerto, I.
Estimation of the offspring mean in a controlled branching process with a random control function. {\it Stochastic Processes and their
Applications}, 117: 928-946, 2007.

\bibitem{ww} Wei, C. Z. and Winnicki, J. Some asymptotic results for the branching process
with immigration. {\it Stochastic Processes and their Applications}, 31(2): 261–282, 1989.

\end{thebibliography}
\end{document}